\documentclass[fleqn]{mat01}
\usepackage{times,mathtimy,amssymb,latexsym,amsbsy}
\begin{document}

\setcounter{page}{217}
\firstpage{217}

\newcommand{\CC}{\mathbb C}
\newcommand{\DD}{\mathbb D}
\newcommand{\HH}{\mathbb H}
\newcommand{\PP}{\mathbb P}
\newcommand{\TT}{\mathbb T}
\newcommand{\ZZ}{\mathbb Z}
\newcommand{\EE}{\mathbb E}
\newcommand{\Znu}{\ZZ^\nu}
\newcommand{\Znuplus}{\ZZ^{\nu +1}_+}
\newcommand{\cc}{\mathcal C}
\newcommand{\dd}{\mathcal D}
\newcommand{\hh}{\mathcal H}
\newcommand{\ww}{\mathcal W}
\newcommand{\WW}{\mathcal S}

\newtheorem{theore}{Theorem}
\renewcommand\thetheore{\arabic{section}.\arabic{theore}}
\newtheorem{theor}[theore]{\bf Theorem}
\newtheorem{rem}[theore]{Remark}
\newtheorem{propo}[theore]{\rm PROPOSITION}
\newtheorem{lem}[theore]{Lemma}
\newtheorem{definit}[theore]{\rm DEFINITION}
\newtheorem{assump}[theore]{Assumption}
\newtheorem{coro}[theore]{\rm COROLLARY}
\newtheorem{exampl}[theore]{Example}
\newtheorem{pot}[theore]{Proof of Theorem}
\newtheorem{case}{Case}

\def\corol{\trivlist \item[\hskip \labelsep{COROLLARY.}]}
\def\noteproof{\trivlist \item[\hskip \labelsep{\it Note added in Proof.}]}

\renewcommand{\theequation}{\thesection\arabic{equation}}

\def\thet{{\bf \vartheta}}


\newcommand{\RR}{\mathbf R}
\newcommand{\rr}{\mathbf{R}}
\newcommand{\eps}{\varepsilon}
\newcommand{\jap}[1]{\langle{#1}\rangle}

\title{New criteria to identify spectrum}

\markboth{A Jensen and M Krishna}{New criteria to identify spectrum}

\author{A JENSEN and M KRISHNA$^{*}$}

\address{Department of Mathematical Sciences and
MaPhySto\thanks{MaPhySto~--~A Network in Mathematical
Physics and Stochastics, funded by the Danish National Research Foundation.}\ \,,
Aalborg University,\\ \noindent Fr. Bajers Vej 7G, DK-9220 Aalborg \O{}, Denmark \\
\noindent $^{*}$Institute of Mathematical Sciences, Taramani, Chennai~600~113, India\\
\noindent E-mail: krishna@imsc.res.in}

\volume{115}

\mon{May}

\parts{2}

\pubyear{2005}

\Date{MS received 14 July 2004; revised 16 March 2005}

\begin{abstract}
In this paper we give some new criteria for identifying the components
of a probability measure, in its Lebesgue decomposition. This enables us
to give new criteria to identify spectral types of self-adjoint
operators on Hilbert spaces, especially those of interest.
\end{abstract}


\maketitle

\vspace{-1pc}

\section{Introduction}

Let us briefly motivate our interest in determining the spectral
type of a self-adjoint operator.

Let $\mu$ be a probability measure on the real line $\mathbf{R}$. It is
well-known that this measure has a Lebesgue decomposition $\mu=\mu_{\rm
ac}+\mu_{\rm sc}+\mu_{\rm pp}$, where $\mu_{\rm ac}$ is absolutely
continuous with respect to the Lebesgue measure on $\mathbf{R}$,
$\mu_{\rm sc}$ is singular with respect to Lebesgue measure, and has no
atomic part, i.e. $\mu_{\rm sc}(\{x\})=0$ for all $x\in\mathbf{R}$, and
$\mu_{\rm pp}$ is purely atomic.

This decomposition of a probability measure has important applications
in the theory of a self-adjoint operator $H$ on a (separable) Hilbert
space $\mathcal{H}$. Associated with $H$ is the spectral measure
$E(\cdot)$. The spectral theorem states that we have
\begin{equation*}
\langle u,Hu\rangle=\int_{\mathbf{R}}\lambda
d\langle u,E(\lambda)u\rangle.
\end{equation*}
If $\|u\|=1$, then $d\langle u,E(\cdot)u\rangle$ is a probability
measure, which is supported on the spectrum $\sigma(H)$ of $H$. The
Lebesgue decomposition of probability measures leads to an orthogonal
decomposition of the Hilbert space
\begin{equation*}
\mathcal{H}=\mathcal{H}_{\rm ac}\oplus\mathcal{H}_{\rm sc}\oplus
\mathcal{H}_{\rm pp}.
\end{equation*}
Each subspace is the closure of vectors $u$, such that $d\langle
u,E(\cdot)u\rangle$ is purely absolutely continuous, etc. The subspaces
reduce the operator $H$, such that $H_{|_{\mathcal{H}_{\rm ac}}}$ is a
self-adjoint operator on $\mathcal{H}_{\rm ac}$, etc. In the case of
absolutely continuous and singular continuous parts, one defines the
corresponding parts of the spectrum to be those of the restricted
operators. In the point spectrum case one usually takes $\sigma_{\rm
pp}(H)$ to be the set of eigenvalues of $H$, in order to handle the
case, when the operator has a dense set of eigenvalues. The spectrum of
the operator restricted to $\mathcal{H}_{\rm pp}$ is then the closure of
this set.\newpage

The spectral types of an operator $H$, which is the Hamiltonian of a
quantum mechanical system, is related to the dynamics of the system,
although the relation is by no means simple. The relation comes from the
representation of the time evolution operator\break $\hbox{e}^{-itH}$ as
\begin{equation*}
\langle u,\hbox{e}^{-itH}u\rangle=\int_{\mathbf{R}}\hbox{e}^{-it\lambda}
d\langle u,E(\lambda)u\rangle.
\end{equation*}
In some quantum mechanical systems (e.g. atoms and molecules) the
absolutely continuous part is related to the scattering states, since
$\langle u,\hbox{e}^{-itH}u\rangle$ tends to zero for $u\in\mathcal{H}_{\rm
ac}$ (a consequence of the Riemann--Lebesgue lemma), and the eigenvalues
of $H$ are related to the bound states. In many of these systems one
expects that the singular continuous component is absent, and many
techniques have been developed to prove this type of result. In solid
state physics the situation is somewhat different, and here one has a
wider variety of spectral types.

These applications have motivated us to seek new criteria for
identifying the spectral type of a self-adjoint operator.

The components of a probability measure can be identified via a
transform of the measure. Two of these are well-known, viz. the Fourier
transform and the Borel transform. In this paper we address the question
of identifying the components using a more general transform. We give
results using a general approximate identity, and an associated
continuous wavelet transform.

Concerning the literature, the connection between an approximate
identity and the continuous wavelet transform was discussed by
Holschneider~\cite{hol95}, while wavelet coefficients of fractal
measures were studied by Strichartz \cite{MR93c:42036}. In the theory
of self-adjoint operators finer decomposition of spectra with respect to
Hausdorff measures was first used by Last \cite{yl} and general criteria
for recovering a measure from its Borel transform was done by Simon
\cite{si1}.

\section{The criteria}

We need to introduce conditions on our function $\psi$. Several of these
can be relaxed in some of the results. We use the standard notation
$\jap{x}=(1+x^2)^{1/2}$.

\begin{assump}\label{a1.1}
{\rm Assume that $\psi\in C^1(\rr)$, $\psi(0)=1$, $\psi$ is even, and there exist
$C>0$ and $\delta>1$, such that
\begin{equation}\label{eqassume}
|\psi(x)|+|x\psi'(x)|\leq C \jap{x}^{-\delta},\quad x\in\rr.
\end{equation}
We set $A_\psi = \int_\rr \psi(x) \hbox{d}x$ and assume that $A_\psi \neq 0$.}
\end{assump}
In the sequel we always impose this assumption on $\psi$. We introduce
the notation
\begin{equation}
\psi_a(x)=\psi(x/a)\quad\hbox{and}\quad \tilde{\psi}_a(x)=\frac{1}{a}\psi_a(x),
\quad a>0.
\end{equation}
In particular, the family $\{A_{\psi}^{-1}\tilde{\psi}_a\}$ is an
approximate identity. Let $\mu$ be a probability measure on $\rr$ in
what follows, with Lebesgue decomposition $\mu=\mu_{\rm s}+\mu_{\rm
ac}$. Let $f$ be a function. We recall that the convolution $(f\ast
\mu)(x)=\int f(x-y)\hbox{d}\mu(y)$ is defined, when the integral converges.
Since $\psi$ is bounded, the convolution $\psi_a\ast\mu$ is defined for
all $a>0$.

For $0\leq\alpha\leq1$ we define
\begin{equation}
(d_{\alpha}\mu)(x)=\lim_{\eps\downarrow0}
\frac{\mu((x-\eps,x+\eps))}{(2\eps)^{\alpha}},
\end{equation}
whenever the limit on the right-hand side exists.

We can now state the results. We first give results based on $\psi_a$
and $\tilde{\psi}_a$, and then on an associated continuous wavelet
transform.

\begin{theor}[\!]
\label{thm1}
Let $\mu$ be a probability measure. Then we have as follows{\rm :}

\begin{enumerate}
\renewcommand\labelenumi{\rm \arabic{enumi}.}
\leftskip -.2pc
\item Let $\psi$ satisfy Assumption~{\rm\ref{a1.1}}. Then for every
continuous function $f$ of compact support{\rm ,} the following is
valid.
\begin{equation*}
\hskip -1.2pc \lim_{a \rightarrow 0} \int (\tilde{\psi}_a * \mu)(x) f(x)
\hbox{\rm d}x = A_\psi \int f(x) \hbox{\rm d}\mu(x).
\end{equation*}
\item $\displaystyle{\lim_{a \rightarrow 0}}
(\psi_{a} * \mu)(x) = \mu(\{x\})$.

\item Assume $0<\alpha\leq1$ and $(d_\alpha \mu)(x)$ finite. Then we
have
\begin{equation}
\label{eqn1.1}
\hskip -1.2pc \lim_{a \rightarrow 0} a^{-\alpha} (\psi_{a}*\mu)(x) =
c_\alpha (d_\alpha\mu)(x),
\end{equation}
where
$c_\alpha = \int_0^{\infty} \alpha 2^\alpha y^{\alpha-1}\psi(y)\hbox{\rm d}y$.
\end{enumerate}
\end{theor}

\begin{rem}$\left.\right.$\vspace{-.5pc}
{\rm \begin{enumerate}
\renewcommand\labelenumi{(\arabic{enumi})}
\leftskip .15pc
\item Equation \eqref{eqn1.1} implies that if $\mu$ is purely singular,
then the limit of $\tilde{\psi}_a * \mu(x)$ is zero almost everywhere
with respect to the Lebesgue measure, since the derivative $(d_1\mu)(x)
= 0$ almost everywhere for purely singular $\mu$.

\item If $x$ is not in the topological support of $\mu$, then for each
$0\leq \alpha \leq 1$,
\begin{equation*}
\hskip -1.2pc \lim_{a \rightarrow 0} a^{-\alpha} \psi_{a}*\mu(x) = 0.
\end{equation*}

Our next theorem says a bit more and the first part is analogous to
Wiener's theorem and its extension by Simon \cite{si1}.
\end{enumerate}}
\end{rem}

\begin{theor}[\!]
\label{thm2}
Let $\mu$ be a probability measure. Then for any bounded interval
$(c,d)$ the following are valid.
\begin{enumerate}
\renewcommand\labelenumi{\rm \arabic{enumi}.}
\leftskip -.2pc
\item Let
\begin{equation*}
\hskip -1.2pc C = \int_\rr |\psi(x)|^2\hbox{\rm d}x,
\end{equation*}
then
\begin{align}
&\hskip -1.2pc \lim_{a \rightarrow 0} \frac{1}{a} \int_c^d |(\psi_a * \mu)(x)|^2 \hbox{\rm d}x\nonumber\\[.3pc]
&\hskip -1.2pc \quad\, = C \left(\sum_{x \in (c,d)} \mu(\{x\})^2 +
\frac{1}{2} [\mu(\{c\})^2 + \mu(\{d\})^2 ]\right).
\end{align}

\item For $0 < p < 1${\rm ,} we have
\begin{equation}
\hskip -1.2pc \lim_{a \rightarrow 0} \int_c^d |(\tilde{\psi}_a * \mu) (x)|^p \hbox{\rm d}x
= |A_{\psi}|^p\int_c^d \left|\frac{{\rm d}\mu_{\rm ac}}{{\rm d}x}(x)\right|^p \hbox{\rm d}x.
\end{equation}
\end{enumerate}
\end{theor}

This theorem has the following corollary.

\begin{coro}$\left.\right.$\vspace{.5pc}
\label{cor1}

\noindent Let $\mu$ be a probability measure. Then we have the following
results{\rm :}
\begin{enumerate}
\renewcommand\labelenumi{\rm \arabic{enumi}.}
\leftskip -.2pc
\item $\mu$ has no point part in $[c,d]${\rm ,} if and only if
\begin{equation}
\hskip -1.2pc \liminf_{a\to0}\frac{1}{a}\int_c^d|(\psi_a\ast\mu)(x)|^2 \hbox{\rm d}x=0.
\end{equation}
\item If $\mu$ has no absolutely continuous part in $(c,d)${\rm ,} if
and only if for some $p, 0<p<1${\rm ,}
\begin{equation}
\hskip -1.2pc \liminf_{a\to0}\int_c^d|(\tilde{\psi}_a\ast\mu)(x)|^p
\hbox{\rm d}x=0.
\end{equation}
\end{enumerate}\vspace{-1pc}
\end{coro}

Now to state the results in terms of the continuous wavelet transform,
we introduce
\begin{equation}\label{cwt}
h(x)=\psi(x)+x\psi'(x).
\end{equation}
Under Assumption~\ref{a1.1} we clearly have
\begin{equation}\label{hdecay}
|h(x)|\leq C\jap{x}^{-\delta},
\end{equation}
with the $\delta$ from the assumption. Integration by parts and
eq.~\eqref{cwt} imply that $h$ satisfies the admissibility condition for a
continuous wavelet, i.e. $\int_{-\infty}^{\infty}h(x)\hbox{d}x=0$.

Thus we can define the continuous wavelet transform of a probability
measure $\mu$ as
\begin{equation}\label{cwttrans}
W_h(\mu)(b,a)=\frac{1}{a}\int_{-\infty}^{\infty}h((b-y)/a)\hbox{d}\mu(y).
\end{equation}
The connection between the approximate identity and this transform is
\begin{equation}\label{relation}
-a\frac{\partial}{\partial a}(\tilde{\psi}_a*\mu)(b)=W_h(\mu)(b,a).
\end{equation}
This result follows from
\begin{equation*}
-a\frac{\partial}{\partial a}\left(\frac{1}{a}\psi\left(\frac{x}{a}\right)\right)=
\frac{1}{a}\left(\psi\left(\frac{x}{a}\right)+\frac{x}{a}\psi'\left(\frac{x}{a}\right)\right),
\end{equation*}
and the definitions.

We have the following analogue of Theorem~\ref{thm1}.

\begin{theor}[\!]\label{thm3}
Let $\mu$ be a probability measure. Then we have the following results{\rm :}
\begin{enumerate}
\renewcommand\labelenumi{\rm \arabic{enumi}.}
\leftskip -.2pc
\item We have
\begin{equation}
\hskip -1.2pc \lim_{\varepsilon\downarrow0}\varepsilon
\int_{\varepsilon}^{\infty}W_h(\mu)(b,a)\frac{{\rm d}a}{a}
=\mu(\{b\}).
\end{equation}
\item Let $0<\alpha\leq 1$. Assume that $(d_{\alpha}\mu)(b)$ exists. Then
\begin{equation}
\hskip -1.2pc \lim_{\varepsilon\downarrow0}\varepsilon^{1-\alpha}
\int_{\varepsilon}^{\infty}W_h(\mu)(b,a)\frac{{\rm d}a}{a}
=c_{\alpha}(d_{\alpha}\mu)(b),
\end{equation}
where $c_{\alpha}$ was defined in Theorem~{\rm\ref{thm1}}.
\end{enumerate}
\end{theor}

\begin{rem}
{\rm We note that for $0<\alpha<1$ we can replace $\int_{\eps}^{\infty}$ by
$\int_{\eps}^M$ for any $M>0$ (see the proof of the Theorem).}
\end{rem}

We also have the following analogue of Theorem~\ref{thm2}(1).

\begin{theor}[\!]\label{thm4}
Let $\mu$ be a probability measure. Then for any bounded interval
$(c,d)$ we have the following result. Let
\begin{equation*}
C_h= \int_{\rr}|h(x)|^2 \hbox{\rm d}x.
\end{equation*}
Then we have
\begin{align}
&\lim_{a\downarrow 0}\int_c^d|W_h(\mu)(b,a)|^2 \hbox{\rm d}b\nonumber\\[.3pc]
&\quad\,=
C_h\left(\sum_{x\in(c,d)}\mu(\{x\})^2 +\frac{1}{2}(\mu(\{c\})^2 +
\mu(\{d\})^2)\right).
\end{align}
\end{theor}

Even when the quantity $(d_\alpha \mu)(x)$ does not exist, it is
possible to say something on the wavelet transforms, to cover the cases
of measures which are not supported on the sets where such limits exist.
Set
\begin{align*}
C_{\mu, \psi}^\alpha(x) =
\displaystyle\limsup_{a \rightarrow 0}
\frac{\psi_a*\mu}{a^\alpha}(x) \quad \hbox{and} \quad
D_\mu^\alpha(x) = \displaystyle\limsup_{\varepsilon \rightarrow 0}
\frac{\mu((x-\varepsilon, x+\varepsilon))}{(2\varepsilon)^\alpha}.
\end{align*}
Then we have the following theorem.

\begin{theor}[\!]
\label{thm5}
Let $\mu$ be a probability measure{\rm ,} and let $\psi$ satisfy
Assumption~{\rm\ref{a1.1}}. Then $C_{\mu, \psi}^\alpha(x)$ is finite for
any $x${\rm ,} whenever $D_\mu^\alpha(x)$ is finite for the same $x${\rm
,} and{\rm ,} if $\psi$ is non-negative{\rm ,} they are both finite or
both infinite.
\end{theor}

\begin{rem}
\label{rem-new-1}
{\rm The above theorem implies that if $\limsup_{a \to0}
|(\tilde{\psi}_a*\mu)(x)| < \infty $ for all $x \in (c, d)$, then there
is no singular part of $\mu$ supported in $(c, d)$.}
\end{rem}

Finally as an application of the above theorems we consider $\hh$ to be
a separable Hilbert space and $A$ a self-adjoint operator. Then we have
the following theorem.

\begin{theor}[\!]
\label{thm6}
Suppose $A$ is a self-adjoint operator on $\hh$. Consider a function
$\psi$ satisfying Assumption~{\rm\ref{a1.1}}. Then
\begin{enumerate}
\renewcommand\labelenumi{\rm \arabic{enumi}.}
\leftskip -.2pc
\item $\lambda$ is in the point spectrum of $A${\rm ,} if for some $f \in
\hh, \|f\| = 1${\rm ,}
\begin{equation*}
\hskip -1.2pc \lim_{a\rightarrow 0} \jap{f, \psi_a(A-\lambda)f} = 0.
\end{equation*}
\item Let $B \subset \rr$ be a Borel set of positive Lebesgue measure.
Then $B \cap \sigma_{\rm ac}(A) \neq \emptyset${\rm ,} if for some $f \in
\hh${\rm ,} $ \|f\|=1${\rm ,}
\begin{equation*}
\hskip -1.2pc \lim_{a \rightarrow 0} \jap{f, \tilde{\psi}_a(A-\lambda)f}
\neq 0, \quad\hbox{\rm for a.e.}\ \lambda\in B.
\end{equation*}
\item The point spectrum of $A$ in $(c, d)$ is empty{\rm ,} if and only
if for some orthonormal basis $\{f_n\}${\rm ,} of $\hh${\rm ,} one has
for every $n${\rm ,}
\begin{equation*}
\hskip -1.2pc  \liminf_{a \rightarrow 0} \frac{1}{a}
\int_c^d |\jap{f_n, \psi_a(A-\lambda)f_n}|^2  \hbox{\rm d}\lambda = 0.
\end{equation*}
\item The absolutely continuous spectrum of $A$ in $(c, d)$ is empty{\rm
,} if and only if for some orthonormal basis $\{f_n\}$ of $\hh${\rm ,}
one has for every $n$ and some $0 < p < 1${\rm ,}
\begin{equation*}
\hskip -1.2pc  \liminf_{a \rightarrow 0}
\int_c^d \left|\frac{1}{a}\jap{f_n, \psi_a(A-\lambda)f_n}\right|^p  \hbox{\rm d}\lambda = 0.
\end{equation*}
\end{enumerate}
\end{theor}

\section{Proofs}

Throughout the computations below the letter $C$ denote a constant,
whose value may vary from line to line.

\setcounter{section}{2}
\setcounter{theore}{1}
\begin{pot}$\left.\right.$\vspace{.5pc}

\noindent {\rm
Part (1): Since $f$ is a continuous function of compact support and $\psi_a$
is bounded for each $a >0$, $f(x) \psi_a(x-y)$ is absolutely integrable
and the integral is uniformly bounded in $y \in \rr$. Therefore, by an
application of Fubini, a change of variable $x \rightarrow ax+y$ and
dominated convergence theorem, in that order, it follows that
\begin{align*}
\lim_{a\rightarrow 0 } \int \hbox{d}x f(x) (\tilde{\psi}_a*\mu)(x)
&=
\lim_{a\rightarrow 0 } \int \hbox{d}x f(x) \int \tilde{\psi}_a(x-y) \hbox{d}\mu(y)\\[.3pc]
&= \lim_{a\rightarrow 0} \int \hbox{d}\mu(y) \int f(x) \tilde{\psi}_a(x-y) \hbox{d}x \\[.3pc]
&= \lim_{a\rightarrow 0} \int \hbox{d}\mu(y) \int f(ax+y) \psi(x) \hbox{d}x \\[.3pc]
&= \int \hbox{d}\mu(y) \int \left(\lim_{a\rightarrow 0} f(ax+y)\right) \psi(x) \hbox{d}x \\[.3pc]
&=\int f(y)\hbox{d}\mu(y) \cdot \int \psi(x) \hbox{d}x.
\end{align*}
Part (2): This is a direct consequence of the definition of the integral
noting pointwise that we have
\begin{equation*}
\lim_{a \rightarrow 0} \psi_{a}(x) =
\begin{cases}
0, &\hbox{if}\ x \neq 0, \\
1, &\hbox{if}\ x = 0.
\end{cases}
\end{equation*}
We also need to use the dominated convergence theorem to interchange the
limit and the integral.

\noindent Part (3): Let $\Phi_\mu$ denote the distribution function of $\mu$. Then we
have
\setcounter{section}{3}
\setcounter{equation}{0}
\begin{align}
\frac{1}{a^\alpha}\int_\rr \psi_a(x-y) \hbox{d}\mu(y)
&= - \frac{1}{a^\alpha}\int_\rr \frac{\rm d}{{\rm d}y}\psi((x-y)/a)
\Phi_\mu(y) \hbox{d}y \nonumber\\[.3pc]
&= \frac{1}{a^\alpha}\int_\rr \psi^\prime (y) \Phi_\mu(x-ay) \hbox{d}y\nonumber\\[.3pc]
&= - \int_0^\infty \!\!\psi^\prime (y)
(2y)^\alpha \frac{\Phi_\mu(x+ay) - \Phi_\mu(x - ay)}{(2ay)^\alpha} \hbox{d}y,
\label{eqnnt}
\end{align}
where in the first step we used integration by parts, in the next step
we used changed variables and in the last step we used the oddness of
$\psi^\prime$ to split the integral into positive and negative half-lines and
multiplied by appropriate powers.

\setcounter{section}{2}
We observe that
\begin{equation*}
(d_\alpha \mu)(x) =  {\lim_{a \rightarrow 0}}
\frac{\Phi_\mu(x + ay) - \Phi_\mu(x- ay)}{(2ay)^\alpha}
\end{equation*}
for each $y \in \rr$, and is finite by assumption. Furthermore, the
function $(\Phi_\mu(x + ay) - \Phi_\mu(x- ay))(2ay)^{-\alpha} $ is a
bounded measurable function, and due to \eqref{eqassume} we can
take the limits inside the integral sign in \eqref{eqnnt} and use the
dominated convergence theorem.

Now doing an integration by parts gives the value of the integral as
stated in the theorem.}
\end{pot}

\setcounter{theore}{3}
\begin{pot}$\left.\right.$\vspace{.5pc}

\noindent {\rm Part (1): We have
\begin{align*}
\frac{1}{a}\int_c^d |\psi_a * \mu(x)|^2 \hbox{d}x
= \iint \!\hbox{d}\mu(y_1) \hbox{d}\mu(y_2) \int_c^d \hbox{d}x
\frac{1}{a}\overline{\psi_a(x - y_1)} \psi_a(x-y_2).
\end{align*}
Since the function $\psi_a$ is bounded, the interval $(c,d)$ is bounded,
and $\mu$ is a probability measure, the right-hand side integral
converges absolutely, so we used Fubini to interchange integrals to get
the equality above. Let
\begin{equation*}
h_a(y_1, y_2) = \int_c^d \hbox{d}x \frac{1}{a}\overline{\psi_a(x - y_1)} \psi_a(x-y_2).
\end{equation*}
Suppose $y_1 \neq y_2$, then using the bound
$|\psi(x)| \leq C\jap{x}^{-\delta}$, we see that the bound
\begin{align*}
|h_a(y_1, y_2)| &\leq \frac{C}{a}\int_{-\infty}^{\infty}
\jap{(x+y_2-y_1)/a}^{-\delta}
\jap{x/a}^{-\delta}\hbox{d}x\\[.3pc]
&= \frac{C}{a}
\left(\int_{|x|\leq|y_1-y_2|/2}+\int_{|x|\geq|y_1-y_2|/2}\right)(\cdots)\hbox{d}x\\[.3pc]
&\leq
\frac{Ca^{\delta}}{|y_1-y_2|^{\delta}}
\int_{-\infty}^{\infty}\jap{x/a}^{-\delta}
\hbox{d}(x/a)\\[.3pc]
&\leq \frac{Ca^{\delta}}{|y_1-y_2|^{\delta}}
\end{align*}
is valid. It follows that $\lim_{a\to0}h_a(y_1,y_2)=0$ for $y_1\neq
y_2$. It remains to consider $y_1=y_2$. This is done by noting that
\begin{equation*}
h_a(y_1, y_1) = \int_c^d   \frac{1}{a} |\psi_a(x - y_1)|^2 \hbox{d}x
 = \int_{(c-y_1)/a}^{(d-y_1)/a}   |\psi(x)|^2 \hbox{d}x,
\end{equation*}
from which taking limits, we obtain the stated value for the
coefficient, either $C$ or $C/2$, based on whether $c < y_1 < d$ or $y_1
= c, d$, using the evenness of $\psi$. Now to complete the proof, we
note the estimate
\begin{equation*}
|h_a(y_1, y_2)|\leq C\int_{\rr}\jap{x/a}^{-\delta}\hbox{d}(x/a)\leq C_0,
\end{equation*}
where the constant $C_0$ is independent of $a$, $y_1$, and $y_2$. Thus
the proof is completed using the dominated convergence theorem.

\noindent Part (2): We adapt the arguments in \cite{si1} to the case at hand. We split
the measure in three components: $\mu=\mu_1+\mu_2+\mu_3$. Here
$\hbox{d}\mu_1=(1-\chi_{[c-1,d+1]})\hbox{d}\mu$, $\hbox{d}\mu_2=g\hbox{d}x$
with $g\in L^1([c-1,d+1])$, and $\mu_3$ is purely singular, and
supported on $[c-1,d+1]$. We have for $x\in[c,d]$ the estimate
\begin{equation*}
|(\tilde{\psi}_a\ast\mu_1)(x)|\leq
C\int_{\rr\setminus[c-1,d+1]}a^{-1}\jap{(x-y)/a}^{-\delta}\hbox{d}\mu_1(y)
\leq C a^{\delta-1}.
\end{equation*}
We now look at the $\mu_2$ part. We have, for $0 < p <1$, by the reverse
H\"older inequality
\begin{align*}
&\int_c^d |(\tilde{\psi}_a *g)(x) - A_\psi g(x)|^p \hbox{d}x\\[.3pc]
&\quad\, \leq
\left(\int_c^d |
(\tilde{\psi}_a *g)(x) - A_\psi g(x)| \hbox{d}x\right)^{p} (d-c)^{1-p},
\end{align*}
which implies that $\tilde{\psi}_a*g \rightarrow A_{\psi} g$ in
$L^p((c, d))$, $0 < p \leq 1$.}
\end{pot}

Now we will show that the singular part $\mu_3$ does not contribute to
the limit. So assume that $\mu_3$ is purely singular and that its
support $S$ is contained in $[c-1, d+1]$. Since $\mu_3$ is singular, by
the definition of support, $S$ satisfies $\mu_3(\RR\hbox{$\setminus$}S) = 0$
and $|S| = 0$, with $|\cdot|$ denoting the Lebesgue measure. By the
regularity of the Lebesgue measure, given an $\varepsilon >0$, there is
an open set $O\subset(c-2,d+2)$, such that $S \subset O$, with
$|O\hbox{$\setminus$}S| < \varepsilon$. We also have $|O| \leq |O\hbox{$\setminus$}S| +
|S| < \varepsilon$. For the same $\varepsilon$, since the measure
$\mu_3$ is regular, we also have a compact $K \subset S$, such that
$\mu_3(S\hbox{$\setminus$}K) < \varepsilon.$ In addition, since $K \subset S$,
and $S$ has Lebesgue measure zero, $K$ also has Lebesgue measure zero.

The above reverse H\"older inequality gives
\begin{align*}
\int_c^d |(\tilde{\psi}_a*\mu_3)(x)|^p \hbox{d}x
&= \int_O |(\tilde{\psi}_a*\mu_3)(x)|^p \hbox{d}x\\[.3pc]
&\quad\,+\int_{(c,d)\setminus O} |(\tilde{\psi}_a*\mu_3)(x)|^p \hbox{d}x\\[.3pc]
&\leq |O|^{1-p} \mu_3((c,d))^p \|\psi\|_1^p \\[.3pc]
&\quad\,+|d-c|^{1-p} \left(\int_{(c,d)\setminus O}
|(\tilde{\psi}_a*\mu_3)(x)| \hbox{d}x\right)^p 
\end{align*}
\begin{align*}
&\leq C \varepsilon^{1-p} + |d-c|^{1-p}\\[.3pc]
&\quad\,\times \left(\int_{(c,d)\setminus O}
|(\tilde{\psi}_a*\mu_3)(x)| \hbox{d}x\right)^p.
\end{align*}
Now consider a bounded continuous function $h$ which is $1$ on $(c,d)
\hbox{$\setminus$}O$, and $0$ on $K$.

Then using Assumption~\ref{a1.1} that $|\psi(x)| \leq C
\jap{x}^{-\delta}$, and setting $\phi(x) = \jap{x}^{-\delta}$,
\begin{align*}
\int_{(c,d)\setminus O} |(\tilde{\psi}_a*\mu_3)(x)| \hbox{d}x
&\leq \int_{(c,d) \setminus O} \frac{1}{a} \int_\RR |\psi_a(x-y)| \hbox{d}\mu_3(y) \hbox{d}x\\[.3pc]
&\leq C\int_{(c,d) \setminus O} \frac{1}{a}
\int_\RR \langle{(x-y)/a}\rangle^{-\delta} \hbox{d}\mu_3(y) \hbox{d}x\\[.3pc]
&\leq C\int_{(c,d) \setminus O}  h(x)(\tilde{\phi}_a*\mu_3)(x) \hbox{d}x.
\end{align*}
The function $\phi$ satisfies Assumption~\ref{a1.1}, so Theorem
\ref{thm1}(1) is applicable with $\psi$ replaced by $\phi$. Therefore
the last term, which has positive integrand, converges to
$\int_{(c,d)\setminus O} h(x) \hbox{d}\mu(x)$ as $a$ goes to zero, which is
bounded by $\int_{(c,d)\setminus K} \hbox{d}\mu(x)$,
\begin{align*}
\int_{(c,d)\setminus O} h(x) d\mu(x) \leq
\mu((c,d)\hbox{$\setminus$}K) \leq
\mu((c,d)\hbox{$\setminus$}S) + \mu(S\hbox{$\setminus$}K) < \varepsilon,
\end{align*}
using the facts that $\mu((c,d)\hbox{$\setminus$}S) = 0$ and $\mu(S\hbox{$\setminus$}K)
< \varepsilon$.

Using the inequality $(a+b+c)^p\leq a^p+b^p+c^p$ for $0<p<1$ and
non-negative numbers $a$, $b$, $c$, we have
\begin{align*}
 \int_c^d |(\tilde{\psi}_a * \mu) (x)-A_{\psi}g(x)|^p  \hbox{d}x
&\leq \int_c^d |(\tilde{\psi}_a * \mu_1) (x)|^p  \hbox{d}x\\[.3pc]
{}&\quad+\int_c^d |(\tilde{\psi}_a * \mu_2) (x)-A_{\psi}g(x)|^p  \hbox{d}x\\[.3pc]
{}&\quad+\int_c^d |(\tilde{\psi}_a * \mu_3) (x)|^p  \hbox{d}x.
\end{align*}
Putting the above estimates together and using $\varepsilon$ arbitrary,
one gets
\begin{equation*}
\lim_{a\rightarrow 0} \int_c^d |(\tilde{\psi}_a * \mu) (x)-A_{\psi}g(x)|^p  \hbox{d}x=0.
\end{equation*}
Now the spaces $L^p((c,d))$, $0<p<1$, are metric spaces with the metric
$d(f,g)=\|f-g\|_p^p$. It then follows from the triangle inequality for
this metric that
\begin{equation*}
\lim_{a \rightarrow 0} \int_c^d |(\tilde{\psi}_a * \mu)(x)|^p \hbox{d}x
= |A_\psi|^p\int_c^d |g(x)|^p \hbox{d}x.
\end{equation*}
Since $g={{\rm d}\mu_{\rm ac}}/{{\rm d}x}$, the result follows.

\setcounter{theore}{5}
\begin{pot}
{\rm Let $0<\eps<M<\infty$. It follows from \eqref{relation} that we have
\begin{equation*}
\int_{\eps}^MW_h(\mu)(b,a)\frac{{\rm d}a}{a}=
(\tilde{\psi}_{\eps}\ast\mu)(b)-(\tilde{\psi}_M\ast\mu)(b).
\end{equation*}
The results now follow from Theorem~\ref{thm1}.}\hfill $\Box$
\end{pot}

\setcounter{theore}{7}
\begin{pot}
{\rm The proof is entirely analogous to the proof of Theorem~\ref{thm2},
replacing $\psi$ by $h$ and adjusting the powers of $a$.} \hfill $\Box$
\end{pot}

\begin{pot}
{\rm Consider the case when $D_\mu^\alpha(x)$ is finite for some $x$ and for
some fixed $\alpha$. Then for any $0 < y < 1$, $\mu(x-y, x+y) \leq C
|y|^\alpha$ for some finite constant $C$. So, using the last line in
eq.~(\ref{eqnnt}) and estimating the right-hand side, one has,
by Assumption~\ref{a1.1},
\begin{align*}
\left|\frac{1}{a^\alpha}(\psi_a*\mu)(x)\right|
\leq C \int_0^{\infty} |\psi^\prime (y)|(2y)^\alpha \hbox{d}y
\leq C \int_0^{\infty} \jap{y}^{-\delta}|y|^{-1+\alpha} \hbox{d}y <\infty.
\end{align*}
Now taking the $\limsup$ of the left-hand side the finiteness of
$C_{\mu, \psi}^\alpha$ follows.

On the other hand, since $\psi$ is positive continuous with $\psi(0) =
1$, there is a $\beta >0 $ such that $\psi(y) > 1/2, -\beta < y <
\beta$. Using this and the evenness~of~$\psi$,
\begin{align*}
\frac{1}{a^\alpha}(\psi_a*\mu)(x)
& = \frac{1}{a^\alpha}\int\! \psi_a(x-y) \hbox{d}\mu(y) = \int \psi(y/a) \hbox{d}\mu(y+x)\\[.3pc]
& \geq \frac{1}{a^\alpha}\int_{-\beta a}^{\beta a} \frac{1}{2} \hbox{d}\mu(y+x) \\[.3pc]
& \geq  \frac{1}{2a^\alpha} [\mu(x + a\beta ) - \mu(x-a\beta)],
\end{align*}
where $\psi\geq 0$ is used to get the first inequality. The above
inequalities immediately imply that since $\beta$ is fixed,
$D_\mu^\alpha(x) =\infty$ implies the same for $C_{\mu,
\psi}^\alpha(x)$.} \hfill $\Box$
\end{pot}

\setcounter{theore}{10}
\begin{pot}
{\rm Parts (1) and (2) are a direct application of Theorems~\ref{thm1}(2) and
(3) respectively. Parts (3) and (4) are a direct application of
Corollaries \ref{cor1}(1) and (2) respectively.}
\end{pot}

\end{document}